\documentclass[12pt]{amsart} 
\usepackage{amssymb,amsmath,amscd} 
\theoremstyle{plain} 
\newtheorem{Lem}{Lemma}[section] 
\newtheorem{Prop}[Lem]{Proposition} 
\newtheorem{Thm}[Lem]{Theorem} 
\newtheorem{Cor}[Lem]{Corollary} 
\theoremstyle{definition} 
\newtheorem{Def}[Lem]{Definition} 
\newtheorem{Rem}[Lem]{Remark} 
 
\newtheorem{Prbl}[Lem]{Problem} 

\errorcontextlines=0  \numberwithin{equation}{section}

\newcommand{\bbC}{{\mathbb C}} 
\newcommand{\bbR}{{\mathbb R}}

\newcommand{\bbN}{{\mathbb N}} 
\newcommand{\bbK}{{\mathbb K}}

\newcommand{\calH}{{\mathcal H}} 
\newcommand{\calL}{{\mathcal L}}

\newcommand{\calP}{{\mathcal P}} 
\newcommand{\calS}{{\mathcal S}} 
\newcommand{\calT}{{\mathcal T}} 
\newcommand{\calV}{{\mathcal V}} 
 
\newcommand{\id}{\mathrm{id}} 
\newcommand{\frakl}{{\mathfrak l}} 
\newcommand{\frakr}{{\mathfrak r}} 

\newcommand{\fR}{{\mathfrak R}}
 
\newcommand{\bpi}{\begin{picture}} 
\newcommand{\epi}{\end{picture}}

\newcommand{\regrep}{\breve\omega}

\begin{document} 

\title[Algebraic covariant derivative curvature tensors] 
{Generators of algebraic covariant derivative curvature tensors and Young symmetrizers} 
\date{September 2003} 
\author[B. Fiedler]{Bernd Fiedler}
\address{Bernd Fiedler \\ Mathematisches Institut \\ Universit\"at Leipzig\\ 
Augustusplatz 10/11 \\ D-04109 Leipzig \\ Germany}
\urladdr{http://home.t-online.de/home/Bernd.Fiedler.RoschStr.Leipzig/}  
\email{Bernd.Fiedler.RoschStr.Leipzig@t-online.de}  

\begin{abstract}
We show that the space of algebraic covariant derivative curvature tensors ${\mathfrak R}'$ is generated by Young symmetrized product tensors $T\otimes \hat{T}$ or $\hat{T}\otimes T$, where $T$ and $\hat{T}$ are covariant tensors of order 2 and 3 whose symmetry classes are irreducible and characterized by the following pairs of partitions: $\{(2),(3)\}$, $\{(2),(2\,1)\}$ or $\{(1^2),(2\,1)\}$. Each of the partitions $(2)$, $(3)$ and $(1^2)$ describes exactly one symmetry class, whereas the partition $(2\,1)$ characterizes an infinite set ${\mathfrak S}$ of irreducible symmetry classes. This set ${\mathfrak S}$ contains exactly one symmetry class $S_0\in{\mathfrak S}$ whose elements $\hat{T}\in S_0$ can not play the role of generators of tensors ${\mathfrak R}'$. The tensors $\hat{T}$ of all other symmetry classes from ${\mathfrak S}\setminus\{S_0\}$ can be used as generators for tensors ${\mathfrak R}'$.

Foundation of our investigations is a theorem of S. A. Fulling, R. C. King, B. G. Wybourne and C. J. Cummins about a Young symmetrizer that generates the symmetry class of algebraic covariant derivative curvature tensors.
Furthermore we apply ideals and idempotents in group rings ${\mathbb C}[{\calS}_r]$, the Littlewood-Richardson rule and discrete Fourier transforms for symmetric groups ${\calS}_r$. For certain symbolic calculations we used the Mathematica packages {\sf Ricci} and {\sf PERMS}.
\end{abstract}

\maketitle 

%
%

\section{Introduction}

In \cite{fie20} we described constructions of generators of {\itshape algebraic curvature tensors}. The present paper searches for generators of {\itshape algebraic covariant derivative curvature tensors}.

Algebraic curvature tensors are tensors of order 4 which have the same symmetry properties as the Riemann tensor of a Levi-Civita connection in Differential Geometry.
Let ${\calT}_r V$ be the vector space of the $r$-times covariant tensors $T$ over a finite-dimensional $\bbK$-vector space $V$, $\bbK = \bbR$ or  $\bbK = \bbC$. We assume that $V$ possesses a {\itshape fundamental tensor} $g \in {\calT}_2 V$ (of arbitrary signature) which can be used for raising and lowering of tensor indices.
\begin{Def}
A tensor $\fR \in {\calT}_4 V$ is called an {\itshape algebraic curvature tensor} iff $\fR$ has the {\itshape index commutation symmetry}
\begin{eqnarray}
\forall\,w, x, y, z\in V:\;\;\;
\fR(w,x,y,z) & = & - \fR(w,x,z,y) \;=\; \fR(y,z,w,x)
\label{equ1.1}
\end{eqnarray}
and fulfills the {\itshape first Bianchi identity}
\begin{eqnarray}
\forall\,w, x, y, z\in V:\;\;\;\fR(w,x,y,z) + \fR(w,y,z,x) + \fR(w,z,x,y) & = & 0\,.
\end{eqnarray}
\end{Def}
\begin{Def}
A tensor $\fR' \in {\calT}_5 V$ is called an {\itshape algebraic covariant derivative curvature tensor} iff $\fR'$ has the {\itshape index commutation symmetry}
\begin{eqnarray}
\fR'(w,x,y,z,u) & = & - \fR'(w,x,z,y,u) \;=\; \fR'(y,z,w,x,u)
\end{eqnarray}
and fulfills the {\itshape first Bianchi identity}
\begin{eqnarray}
\fR'(w,x,y,z,u) + \fR'(w,y,z,x,u) + \fR'(w,z,x,y,u) & = & 0
\end{eqnarray}
and the {\itshape second Bianchi identity}
\begin{eqnarray}
\fR'(w,x,y,z,u) + \fR'(w,x,z,u,y) + \fR'(w,x,u,y,z) & = & 0 \label{equ1.5}
\end{eqnarray}
for all $u, w, x, y, z\in V$.
\end{Def}
The relations (\ref{equ1.1}) -- (\ref{equ1.5}) correspond to the well-known formulas
\begin{eqnarray}
 & & R_{i j k l}\;=\; - R_{i j l k}\;=\;R_{k l i j} \label{equ1.6}\\
 & & R_{i j k l} + R_{i k l j} + R_{i l j k}\;=\;0
\end{eqnarray}
for the Riemann tensor $R$ and
\begin{eqnarray}
 & & R_{i j k l\,;\,m}\;=\; - R_{i j l k\,;\,m}\;=\;R_{k l i j\,;\,m}\\
 & & R_{i j k l\,;\,m} + R_{i k l j\,;\,m} + R_{i l j k\,;\,m}\;=\;0\\
 & & R_{i j k l\,;\,m} + R_{i j l m\,;\,k} + R_{i j m k\,;\,l}\;=\;0
\end{eqnarray}
for its first covariant derivative which we present here in terms of tensor coordinates with respect to an arbitrary local coordinate system.

Investigations of tensors of type $\fR$ and $\fR'$ were carried out by many authors. (See the extensive bibliography in the book \cite{gilkey5} by P. B. Gilkey.) A famous problem connected with algebraic curvature tensors is the {\itshape Osserman conjecture}.
\begin{Def}
Let $\fR \in {\calT}_4 V$ be an algebraic curvature tensor and $x \in V$ be a vector with $|g(x,x)| = 1$.
The {\itshape Jacobi operator} $J_{\fR}(x)$ of $\fR$ and $x$ is the linear operator
$J_{\fR}(x) : V \rightarrow V\;,\; J_{\fR}(x): y \mapsto J_{\fR}(x)y$ that is defined by
\begin{eqnarray}
\forall\, w \in V : \;\;\;
g(J_{\fR}(x) y , w) & = & \fR(y, x, x, w)\,. \label{equ1.11}
\end{eqnarray} 
\end{Def}
%
%
%
\begin{Def}
An algebraic curvature tensor $\fR$ is called {\itshape spacelike Osserman} (resp. {\itshape timelike Osserman}) if the eigenvalues of $J_{\fR}(x)$ are constant on $S^{+}(V):=\{x\in V\,|\,g(x,x)=+1\}$ (resp. $S^{-}(V):=\{x\in V\,|\,g(x,x)=-1\}$).
\end{Def}
Since the notions ''spacelike Osserman'' and ''timelike Osserman'' are equivalent
one simply says {\itshape Osserman}.

If $R$ is the Riemann tensor of a Riemannian manifold $(M,g)$ which is locally a rank one symmetric space or flat, then the eigenvalues of $J_R(x)$ are constant on the unit sphere bundle of $(M,g)$. Osserman \cite{oss90} wondered if the converse held. This question is known as the {\itshape Osserman conjecture}.

The correctness of the Osserman conjecture has been established for Riemannian manifolds $(M,g)$ in all dimensions $\not= 8, 16$ (see \cite{c88,n02}) and for Lorentzian manifolds $(M,g)$ in all dimensions (see \cite{bbg97,gkv97}). However Osserman's question has a negative answer in the case of a pseudo-Riemannian metric with signature $(p,q)$, $p,q\ge 2$ (see e.g. the references given in \cite{fie20}).
A detailed view about the Osserman conjecture can be found in the book \cite{gilkey5} by P. B. Gilkey.

Numerous examples of Osserman algebraic curvature tensors can be constructed by menas of operators $\alpha$ and $\gamma$ given below. It turned out that these operators lead to generators for arbitrary algebraic curvature tensors.
\begin{Def}
\begin{enumerate}
\item{Let $S \in {\calT}_2 V$ be a symmetric tensor of order 2, i.e. the coordinates of $S$ satisfy $S_{i j} = S_{j i}$. We define a tensor $\gamma (S) \in {\calT}_4 V$ by
\begin{eqnarray}
\gamma (S)_{i j k l} & := &
{\textstyle \frac{1}{3}} \left( S_{i l} S_{j k} - S_{i k} S_{j l} \right) \,.
\label{equ1.12}
\end{eqnarray}
}
\item{Let $A \in {\calT}_2 V$ be a skew-symmetric tensor of order 2, i.e. the coordinates of $A$ satisfy $A_{j i} = - A_{i j}$. We define a tensor $\alpha (A) \in {\calT}_4 V$ by
\begin{eqnarray}
\alpha (A)_{i j k l} & := &
{\textstyle \frac{1}{3}} \left( 2\,A_{i j} A_{k l} + A_{i k} A_{j l}
-  A_{i l} A_{j k} \right) \,. \label{equ1.13}
\end{eqnarray}
}
\end{enumerate}
\end{Def}
Now we can construct an example of an Osserman algebraic curvature tensor in the following way. Let $g\in{\calT}_2 V$ be a positive definite metric and $\{ C_i \}_{i = 1}^r$ be a finite set of real, skew-symmetric
$(\mathrm{dim} V \times \mathrm{dim} V)$-matrices that satisfy the Clifford commutation relations
\begin{eqnarray}
C_i \cdot C_j + C_j \cdot C_i & = & - 2\,{\delta}_{i j} \,. \label{equ1.14}
\end{eqnarray}
If we form skew-symmetric tensors $A_i\in\calT_2 V$ by $A_i{(x,y)} := g(C_i\cdot x,y)$ ($x,y\in V$), then
\begin{eqnarray}
\fR & = & {\lambda}_0\,\gamma (g) + \sum_{i = 1}^r
{\lambda}_i \,\alpha (C_i) \;\;\;,\;\;\;{\lambda}_0 , {\lambda}_i = \mathrm{const.}
\end{eqnarray}
is an Osserman algebraic curvature tensor (see \cite{gilkey3}). Further examples which allow also indefinite metrics can be found in \cite[pp.191-193]{gilkey5}. (See also \cite[Sec.6]{fie20}.)

The operators $\alpha$ and $\gamma$ can be used to form generators for arbitrary algebraic curvature tensors. P. Gilkey \cite[pp.41-44]{gilkey5} and B. Fiedler \cite{fie20} gave different proofs for
\begin{Thm}
Each of the sets of tensors
\begin{enumerate}
\item{$\{ \gamma(S)\;|\; S\in\calT_2 V\;\mathrm{symmetric} \}$}
\item{$\{ \alpha(A)\;|\; A\in\calT_2 V\;\mbox{\rm skew-symmetric} \}$}
\end{enumerate}
generate the vector space of all algebraic curvature tensors $\fR$ on $V$.
\end{Thm}
Note that the tensors $\gamma(S)$ and $\alpha(A)$ are expressions which arise from
$S\otimes S$ or $A\otimes A$ by a symmetrization
\begin{eqnarray}
\gamma(S)\;=\;\textstyle{\frac{1}{12}}\,y_t^{\ast}(S\otimes S)
&\;\;\;,\;\;\;&
\alpha(A)\;=\;\textstyle{\frac{1}{12}}\,y_t^{\ast}(A\otimes A)
\end{eqnarray}
where $y_t$ is the Young symmetrizer of the Young tableau
\begin{eqnarray}
t & = 
\begin{array}{|c|c|}
\hline
1 & 3 \\
\hline
2 & 4 \\
\hline
\end{array}
\end{eqnarray}
(see \cite{fie20}).

In the present paper we search for similar {\itshape generators of algebraic covariant derivative curvature tensors}. We use Boerner's definition of {\itshape symmetry classes} for tensors $T\in\calT_r V$ by right ideals $\frakr\subseteq\bbK[\calS_r]$ of the group $\bbK[\calS_r]$ of the symmetric group $\calS_r$ (see Section \ref{sec2} and \cite{boerner,boerner2,fie16,fie18}). On this basis we investigate the following
\begin{Prbl}
We search for generators of algebraic covariant derivative curvature tensors which can be formed by a suitable symmetry operator from tensors
\begin{eqnarray}
T\otimes \hat{T}\;\;\;\mathrm{or}\;\;\;\hat{T}\otimes T &\;\;\;,\;\;\;&
T\in\calT_2 V\;,\;\hat{T}\in\calT_3 V \label{equ1.18}
\end{eqnarray}
where $T$ and $\hat{T}$ belongs to symmetry classes of $\calT_2 V$ and $\calT_3 V$ which are defined by minimal right ideals $\frakr\subset\bbK[\calS_2]$ and $\hat{\frakr}\subset\bbK[\calS_3]$, respectively. \label{probl1.7}
\end{Prbl}
We will see that all such generators can be gained by means of the {\itshape Young symmetrizer} $y_{t'}$ of the {\itshape Young tableau}
\begin{eqnarray}
t' & = &
\begin{array}{|c|c|c|c}
\cline{1-3}
1 & 3 & 5 & \\
\cline{1-3}
2 & 4 &\multicolumn{2}{c}{\;\;\;} \\
\cline{1-2}
\end{array}
\,.
\end{eqnarray}
We obtain the following results:
\begin{Thm} \label{thm1.8}
A solution of Problem {\rm\ref{probl1.7}} can be constructed at most from such pairs of tensors {\rm (\ref{equ1.18})} whose symmetry classes are characterized by the following partitions $\lambda\vdash 2$, $\hat{\lambda}\vdash 3$:
\begin{center}
{\rm
\begin{tabular}{|c|c|c|}
\hline
 & $\lambda$ for $T$ & $\hat{\lambda}$ for $\hat{T}$ \\
\hline
(a) & $(2)$ & $(3)$ \\
(b) & $(2)$ & $(2\,1)$ \\
(c) & $(1^2)$ & $(2\,1)$ \\
\hline
\end{tabular}\,.
}
\end{center}
\end{Thm}
The case (a) of Theorem \ref{thm1.8} is specified by
\begin{Thm}
Let us denote by $S\in \calT_2 V$ and $\hat{S}\in \calT_3 V$ symmetric tensors of order {\rm 2} and {\rm 3}, respectively.
Then the set of all tensors which belong to exactly one of the following tensor types
\begin{eqnarray}
\tau: & &
\begin{array}{cccc}
y_{t'}^{\ast} (S\otimes \hat{S}) & , & y_{t'}^{\ast} (\hat{S}\otimes S) & , \\
\end{array} \label{equ1.20}
\end{eqnarray}
generates the vector space of all algebraic covariant derivative curvature tensors $\fR'\in\calT_5 V$.

Moreover, the tensors {\rm (\ref{equ1.20})} coincide and their coordinates fulfill
\begin{eqnarray} \label{equ1.21}
\hat{\gamma}(S,\hat{S})_{ijkls} & := &
(y_{t'}^{\ast} (S\otimes \hat{S}))_{ijkls} \;=\; (y_{t'}^{\ast} (\hat{S}\otimes S))_{ijkls} \\
 & = & 4\,\left\{S_{il}{\hat{S}}_{jks} - S_{jl}{\hat{S}}_{iks} + S_{jk}{\hat{S}}_{ils} - S_{ik}{\hat{S}}_{jls} \right\} \nonumber
\end{eqnarray} \label{thm1.9}
\end{Thm}
The operator $\hat{\gamma}$ plays the same role for the generators of algebraic covariant derivative curvature tensors considered in Theorem \ref{thm1.9} as the operators $\alpha$ and $\gamma$ play for the generators of algebraic curvature tensors. A first proof that the expressions (\ref{equ1.21}) are generators for $\fR'$ was given by P. B. Gilkey \cite[p.236]{gilkey5}.

The cases (b) and (c) of Theorem \ref{thm1.8} lead to
\begin{Thm} \label{thm1.10}
Let us denote by $S, A \in \calT_2 V$ symmetric or alternating tensors of order {\rm 2} and by $U \in \calT_3 V$ covariant tensors of order {\rm 3} whose symmetry class $\calT_{\frakr}$ is defined by a fixed minimal right ideal $\frakr\subset\bbK[\calS_3]$ from the equivalence class characterized by the partition $(2 , 1) \vdash 3$. We consider the following types $\tau$ of tensors
\begin{eqnarray} \label{equ1.22}%
\tau: & &
\begin{array}{cccc}
y_{t'}^{\ast} (S \otimes U) & , & y_{t'}^{\ast} (U \otimes S) & , \\
y_{t'}^{\ast} (A \otimes U) & , & y_{t'}^{\ast} (U \otimes A) & . \\
\end{array}
\end{eqnarray}
Then for each of the  above types $\tau$ the following assertions are equivalent:
\begin{enumerate}
\item{The vector space of algebraic covariant derivative curvature tensors $\fR' \in \calT_5 V$ is the set of all finite sums of tensors of the type $\tau$ considered. \label{statement1}}
\item{The right ideal $\frakr$ is different from the right ideal
$\frakr_0 := f \cdot \bbK [\calS_3]$ with generating idempotent
\begin{eqnarray}
f \;:=\; \left\{ \frac{1}{2}\,(\id - (1 \,3)) - \frac{1}{6}\,y \right\}
& \;\;\;,\;\;\; &
y \;:=\; \sum_{p \in \calS_3} \mathrm{sign}(p)\,p \,. \label{equ1.23}
\end{eqnarray}
}
\end{enumerate}
\end{Thm}
In the situation of Theorem \ref{thm1.10} we can also determine operators of the type $\alpha$, $\gamma$, $\hat{\gamma}$ which describe the generators of the algebraic covariant derivative curvature tensors $\fR'$ considered. However, these operators depend on the right ideal $\frakr$ (or its generating idempotent $e$) that defines the symmetry class of $U$. And they yield no short expressions of 2, 3, or 4 terms but longer expressions between 10 and 20 terms of length. We describe their determination in a forthcoming paper.

Here is a brief outline to the present paper. In Section \ref{sec2} we give a summary of basic facts about symmetry classes, in particular about the connection between product tensors and Littlewood-Richardson products of corresponding representations. Such Littlewood-Richardson products are used in Subsection \ref{subsec3.1} to prove Theorem \ref{thm1.8}. In the Subsections \ref{subsec3.2} and \ref{subsec3.3} we prove the Theorems \ref{thm1.9} and \ref{thm1.10} by a method which we already applied in \cite{fie20}. Idempotents which define the symmetry class of the tensors $U$ are determined by means of diskrete Fourier transforms.\vspace{20pt}

\section{Basic facts about symmetry classes} \label{sec2}
The vector spaces of algebraic curvature tensors or algebraic covariant derivative tensors over $V$ are {\itshape symmetry classes} in the sence of H. Boerner \cite[p.127]{boerner}. We denote by $\bbK [{\calS}_r]$ the {\itshape group ring} of a symmetric group ${\calS}_r$ over the field $\bbK$. Every group ring element $a = \sum_{p \in {\calS}_r} a(p)\,p \in \bbK [{\calS}_r]$ acts as so-called {\itshape symmetry operator} on tensors $T \in {\calT}_r V$ according to the definition
\begin{eqnarray}
(a T)(v_1 , \ldots , v_r) & := & \sum_{p \in {\calS}_r} a(p)\,
T(v_{p(1)}, \ldots , v_{p(r)}) \;\;\;\;\;,\;\;\;\;\;
v_i \in V \,. \label{equ2.1}
\end{eqnarray}
Equation \eqref{equ2.1} is equivalent to
\begin{eqnarray}
(a T)_{i_1 \ldots i_r} & = & \sum_{p \in {\calS}_r} a(p)\,
T_{i_{p(1)} \ldots  i_{p(r)}} \,.
\end{eqnarray}
\begin{Def}
Let $\frakr \subseteq \bbK [{\calS}_r]$ be a right ideal of $\bbK [{\calS}_r]$ for which an $a \in \frakr$ and a $T \in {\calT}_r V$ exist such that $aT \not= 0$. Then the tensor set
\begin{eqnarray}
{\calT}_{\frakr} & := & \{ a T \;|\; a \in \frakr \;,\;
T \in {\calT}_r V \}
\end{eqnarray}
is called the {\itshape symmetry class} of tensors defined by $\frakr$.
\end{Def}
Since $\bbK [{\calS}_r]$ is semisimple for $\bbK = \bbR , \bbC$, every right ideal $\frakr \subseteq \bbK [{\calS}_r]$ possesses a generating idempotent $e$, i.e. $\frakr$ fulfils $\frakr = e \cdot \bbK [{\calS}_r]$. It holds (see e.g. \cite{fie20} or \cite{boerner,boerner2})
\begin{Lem}
If $e$ is a generating idempotent of $\frakr$, then a tensor $T \in {\calT}_r V$ belongs to ${\calT}_{\frakr}$ iff
\begin{eqnarray}
e T & = & T \,.
\end{eqnarray}
Thus we have
\begin{eqnarray}
{\calT}_{\frakr} & = & \{ eT \;|\; T \in {\calT}_r V \} \,.
\end{eqnarray}
\end{Lem}
Now we summarize tools from our Habilitationsschrift \cite{fie17} (see also its summary \cite{fie18}).
We make use of the following connection
between $r$-times covariant tensors $T \in {\calT}_r V$
and elements of the {\it group ring}
${\Bbb K} [{\calS}_r]$.
\begin{Def} \label{def2.3}
 Any tensor
 $T \in {\calT}_r V$
 and any $r$-tuple
 $b := (v_1 , \ldots , v_r ) \in V^r$
 of
 $r$
 vectors from
 $V$
 induce a function
 $T_b : {\calS}_r \rightarrow {\Bbb K}$
 according to the rule
 \begin{eqnarray}
T_b (p) & := & T(v_{p(1)} , \ldots , v_{p(r)})\;\;\;,\;\;\;p \in {\calS}_r \,.
\end{eqnarray}
We identify this function with the group ring element
$T_b := \sum_{p \in {\calS}_r}T_b (p)\,p \in {\Bbb K} [{\calS}_r]$.
\end{Def}
Obviously,
two tensors $S , T \in {\calT}_r V$ fulfil $S = T$ iff
$S_b = T_b$ for all
$b \in V^r$.
We denote by '$\ast$'
the mapping
$\ast : a = \sum_{p \in {\calS}_r} a(p)\,p \;\mapsto\; a^{\ast} :=
\sum_{p \in {\calS}_r} a(p)\,p^{-1}$. Then the following important formula\footnote{See B. Fiedler \cite[Sec.III.1]{fie16} and B. Fiedler \cite{fie17}.} holds
\begin{eqnarray} \label{equ2.7}
\forall\,T\in\calT_r V\;,\;a\in\bbK[\calS_r]\;,\;b\in V^r\;:\;\;\;\;
(a\,T)_b & = & T_b\cdot a^{\ast}\,.
\end{eqnarray}
Now it can be shown that all $T_b$ of tensors $T$ of a given symmetry class lie in a certain left ideal of ${\Bbb K}[{\calS}_r]$.
\begin{Prop}\hspace{-1mm}\footnote{See B. Fiedler \cite{fie17} or
B. Fiedler \cite[Prop. III.2.5, III.3.1, III.3.4]{fie16}.}
\label{prop2.4}%
Let $e \in {\Bbb K}[{\calS}_r]$ be an idempotent. Then a
$T \in {\calT}_r V$ 
fulfils the condition
$eT = T$
iff
$T_b \;\in\; {\frak l} := {\Bbb K} [{\calS}_r] \cdot e^{\ast}$ for all
$b \in V^r$, i.e.
all $T_b$ of $T$
lie in the left ideal ${\frak l}$ generated by $e^{\ast}$.
\end{Prop}
The proof follows easily from (\ref{equ2.7}). Since a rigth ideal $\frakr$ defining a symmetry class and the left ideal $\frakl$ from Proposition \ref{prop2.4} satisfy $\frakr = \frakl^{\ast}$, we denote symmetry classes also by $\calT_{\frakl^{\ast}}$. A further result is
\begin{Prop}\hspace{-1mm}\footnote{See B. Fiedler \cite{fie17} or
B. Fiedler \cite[Prop. III.2.6]{fie16}.}
\label{prop2.5}%
If $\dim V \ge r$, then every left ideal
${\frak l} \subseteq {\Bbb K}[{\calS}_r]$ fulfils
${\frak l} = {\calL}_{\Bbb K} \{ T_b \;|\;
T \in {\calT}_{{\frak l}^{\ast}} \,,\, b \in V^r \}$.
(Here ${\calL}_{\Bbb K}$ denotes the forming of the linear closure.)
\end{Prop}
If $\dim V < r$, then the $T_b$ of the tensors from
${\calT}_{{\frak l}^{\ast}}$ 
will span only a linear subspace of
${\frak l}$ 
in general.

Important special symmetry operators are Young symmetrizers, which are defined by means of Young tableaux.

A {\itshape Young tableau} $t$ of $r\in\bbN$ is an arrangement of $r$ boxes such that
\begin{enumerate}
\item{the numbers ${\lambda}_i$ of boxes in the rows $i = 1 , \ldots , l$ form a decreasing sequence
${\lambda}_1 \ge {\lambda}_2 \ge \ldots \ge {\lambda}_l > 0$ with
${\lambda}_1 + \ldots + {\lambda}_l = r$,}
\item{the boxes are fulfilled by the numbers $1, 2, \ldots , r$ in any order.}
\end{enumerate}
For instance, the following graphics shows a Young tableau of $r = 16$.
\[\left.
\begin{array}{cc|c|c|c|c|c|c}
\cline{3-7}
{\lambda}_1 = 5 & \;\;\; & 11 & 2 & 5 & 4 & 12 & \\
\cline{3-7}
{\lambda}_2 = 4 & \;\;\; & 9 & 6 & 16 & 15 & \multicolumn{2}{c}{\;\;\;} \\
\cline{3-6}
{\lambda}_3 = 4 & \;\;\; & 8 & 14 & 1 & 7 & \multicolumn{2}{c}{\;\;\;} \\
\cline{3-6}
{\lambda}_4 = 2 & \;\;\; & 13 & 3 & \multicolumn{4}{c}{\hspace{2cm}} \\
\cline{3-4}
{\lambda}_5 = 1 & \;\;\; & 10 & \multicolumn{4}{c}{\hspace{2cm}} \\
\cline{3-3}
\end{array}\right\}\;=\;t\,.
\]
Obviously, the unfilled arrangement of boxes, the {\itshape Young frame}, is characterized by a partition
$\lambda = ({\lambda}_1 , \ldots , {\lambda}_l) \vdash r$ of $r$.

If a Young tableau $t$ of a partition $\lambda \vdash r$ is given, then the {\itshape Young symmetrizer} $y_t$ of $t$ is defined by\footnote{We use the convention $(p \circ q) (i) := p(q(i))$ for the product of two permutations $p, q$.}
\begin{eqnarray}
y_t & := & \sum_{p \in {\calH}_t} \sum_{q \in {\calV}_t} \mathrm{sign}(q)\, p \circ q
\end{eqnarray}
where ${\calH}_t$, ${\calV}_t$ are the groups of the {\itshape horizontal} or
{\itshape vertical permutations} of $t$ which only permute numbers within rows or columns of $t$, respectively. The Young symmetrizers of $\bbK [{\calS}_r]$ are essentially idempotent and define decompositions
\begin{eqnarray}
\bbK [{\calS}_r] \;=\;
\bigoplus_{\lambda \vdash r} \bigoplus_{t \in {\calS\calT}_{\lambda}}
\bbK [{\calS}_r]\cdot y_t
& \;\;,\;\; &
\bbK [{\calS}_r] \;=\;
\bigoplus_{\lambda \vdash r} \bigoplus_{t \in {\calS\calT}_{\lambda}}
y_t \cdot \bbK [{\calS}_r] \label{decomp}
\end{eqnarray}
of $\bbK [{\calS}_r]$ into minimal left or right ideals. In \eqref{decomp}, the symbol ${\calS\calT}_{\lambda}$ denotes the set of all standard tableaux of the partition $\lambda$. Standard tableaux are Young tableaux in which the entries of every row and every column form an increasing number sequence.\footnote{About Young symmetrizers and
Young tableaux see for instance
\cite{boerner,boerner2,full4,fulton,jameskerb,kerber,littlew1,mcdonald,%
muell,naimark,%
waerden,weyl1}. In particular, properties of Young symmetrizers in the case
${\bbK} \not= {\bbC}$ are described in \cite{muell}.}

S.A. Fulling, R.C. King, B.G.Wybourne and C.J. Cummins showed in \cite{full4} that the symmetry classes of the Riemannian curvature tensor $R$ and its {\itshape symmetrized\footnote{$(\,\ldots\,)$ denotes the symmetrization with respect to the indices $s_1, \ldots , s_u$.} covariant derivatives}
\begin{eqnarray}
\left({\nabla}^{(u)}R\right)_{i j k l s_1 \ldots s_u} & := & {\nabla}_{(s_1} {\nabla}_{s_2} \ldots {\nabla}_{s_u)} R_{i j k l}\;=\;R_{i j k l\,;\,(s_1 \ldots s_u)}
\end{eqnarray}
are generated by special Young symmetrizers.
\begin{Thm} \label{thm2.3}
Consider the Levi-Civita connection $\nabla$ of a pseudo-Riemannian metric $g$.
For $u \ge 0$ the Riemann tensor and its symmetrized covariant derivatives
${\nabla}^{(u)} R$ fulfil
\begin{eqnarray}
e_t^{\ast} {\nabla}^{(u)} R & = & {\nabla}^{(u)} R
\end{eqnarray}
where $e_t := y_t (u+1)/(2\cdot (u+3)!)$ is an idempotent which is formed from the Young symmetrizer $y_t$ of the standard tableau
\begin{eqnarray} \label{equ2.12}%
t & = &
\begin{array}{|c|c|c|cc|c|}
\hline
1 & 3 & 5 & \ldots & \ldots & (u+4) \\
\hline
2 & 4 & \multicolumn{4}{l}{\hspace{3cm}} \\
\cline{1-2}
\end{array} \,.
\end{eqnarray}
\end{Thm}
A proof of this result of \cite{full4} can be found in \cite[Sec.6]{fie12}, too. The proof needs only the symmetry properties (\ref{equ1.1}) or (\ref{equ1.6}), the identities Bianchi I and Bianchi II and the symmetry with respect to $s_1, \ldots , s_u$. Thus Theorem \ref{thm2.3} is a statement about algebraic curvature tensors and algebraic covariant derivative curvature tensors.
We can specify this in the following way:
\begin{Def}
A tensor $\fR^{(u)} \in {\calT}_{4+u} V$, $u\ge 0$, is called a {\itshape symmetric algebraic covariant derivative curvature tensor of order $u$} iff $\fR^{(u)}(w,y,z,x,a_1,\ldots,a_u)$ is symmetric with respect to $a_1,\ldots,a_u$ and fulfills
\begin{eqnarray}
\fR^{(u)}(w,x,y,z,a_1,\ldots,a_u) & = & - \fR^{(u)}(w,x,z,y,a_1,\ldots,a_u) \\
 & = & \fR^{(u)}(y,z,w,x,a_1,\ldots,a_u) \nonumber
\end{eqnarray}
\begin{eqnarray}
\hspace*{1cm}0 & = & \fR^{(u)}(w,x,y,z,a_1,\ldots,a_u) + \fR^{(u)}(w,y,z,x,a_1,\ldots,a_u) + \\
 & & \fR^{(u)}(w,z,x,y,a_1,\ldots,a_u) \nonumber \\
0 & = & \fR^{(u)}(w,x,y,z,a_1,a_2,\ldots,a_u) + \fR^{(u)}(w,x,z,a_1,y,a_2,\ldots,a_u) + \\
 & & \fR^{(u)}(w,x,a_1,y,z,a_2,\ldots,a_u) \nonumber
\end{eqnarray}
for all $a_1,\ldots, a_u, w, x, y, z\in V$.
\end{Def}
Now symmetric algebraic covariant derivative curvature tensors can be characterized by means of the Young symmetrizer of the tableau (\ref{equ2.12}).
\begin{Prop} \label{prop2.8a}%
A tensor $T \in {\calT}_{4+u} V$, $u\ge 0$, is a symmetric algebraic covariant derivative curvature tensor of order $u$ iff $T$ satisfies
\begin{eqnarray}
e_t^{\ast} T & = & T
\end{eqnarray}
where $e_t$ is the idempotent from Theorem {\rm \ref{thm2.3}}.
\end{Prop}
A proof of Proposition \ref{prop2.8a} is given in the proof of \cite[Prop.6.1]{fie12}. If we consider now the values $u = 0$ and $u = 1$, we obtain
\begin{Cor} \label{cor2.4}
\begin{enumerate}
\item{A tensor $T \in {\calT}_4 V$ is an algebraic curvature tensor iff $T$ satisfies
\begin{eqnarray}
y_t^{\ast} T & = & 12\,T
\end{eqnarray}
where $y_t$ is the Young symmetrizer of the standard tableau
\begin{eqnarray}
t & = &
\begin{array}{|c|c|}
\hline
1 & 3 \\
\hline
2 & 4 \\
\hline
\end{array} \,.
\end{eqnarray}}
\item{A tensor $\tilde{T} \in {\calT}_5 V$ is an algebraic covariant derivative curvature tensor iff $\tilde{T}$ satisfies
\begin{eqnarray}
y_{t'}^{\ast} \tilde{T} & = & 24\,\tilde{T}
\end{eqnarray}
where $y_{t'}$ is the Young symmetrizer of the standard tableau
\begin{eqnarray} \label{equ2.16}
t' & = &
\begin{array}{|c|c|c|c}
\cline{1-3}
1 & 3 & 5 & \\
\cline{1-3}
2 & 4 &\multicolumn{2}{c}{\;\;\;} \\
\cline{1-2}
\end{array}
\,.
\end{eqnarray}}
\end{enumerate}
\end{Cor}
\vspace{0.5cm}
Now we describe the left ideal $\frakl$ which defines the symmetry class $\calT_{{\frakl}^{\ast}}$ of a product tensor.
We consider products
$T^{(1)} \otimes\ldots\otimes T^{(m)}$ with possibly different $T^{(i)}$.
\begin{Prop}\hspace{- 1mm}\footnote{See B. Fiedler \cite[Sec.III.3.2]{fie16} and B. Fiedler \cite{fie17}.} \label{prop2.8}%
Let ${\frak l}_i \subseteq {\Bbb K} [{\calS}_{r_i}]$ $(i = 1,\ldots , m)$ be left ideals
and
$T^{(i)} \in {\calT}_{{\frak l}_i^{\ast}} \subseteq {\calT}_{r_i} V$ be $r_i$-times 
covariant tensors from the symmetry classes characterized by the
${\frak l}_i$.
Consider the product
\begin{eqnarray}
T & := & T^{(1)} \otimes\ldots\otimes T^{(m)} \;\in\; {\calT}_r V
\;\;\;,\;\;\;
r := r_1 + \ldots + r_m \,. \label{equ2.17}%
\end{eqnarray}
For every $i$ we define an embedding
\begin{eqnarray}
{\iota}_i : {\calS}_{r_i} \rightarrow {\calS}_r
& \;\;,\;\; &
({\iota}_i s)(k) := 
\left\{
\begin{array}{ll}
{\Delta}_i + s(k - {\Delta}_i) & {\rm if}\;\; r_{i-1} < k \le r_i \\
k & {\rm else}
\end{array}
\right.
\end{eqnarray}
where ${\Delta}_i := r_0 + \ldots + r_{i-1}$ and $r_0 := 0$.
Then the $T_b$ of the tensor {\rm (\ref{equ2.17})} fulfil
\begin{eqnarray}
 & \;\;\;\;\; &
\forall\, b \in V^r :\;
 T_b \;\in\; {\frak l} \; := \;
{\Bbb K} [{\calS}_r]\cdot{\calL}\bigl\{{\tilde{\frak l}}_1 \cdot\ldots\cdot {\tilde{\frak l}}_m
\bigr\} \; = \;
{\Bbb K} [{\calS}_r]\cdot\bigl({\tilde{\frak l}}_1 \otimes\ldots\otimes
{\tilde{\frak l}}_m \bigr) \label{equ2.19}%
\end{eqnarray}
where ${\tilde{\frak l}}_i := {\iota}_i ({\frak l}_i)$ are the embeddings of
the ${\frak l}_i$ into ${\Bbb K} [{\calS}_r]$ induced by the ${\iota}_i$. If
$\dim V \ge r$, then the above left ideal ${\frak l}$ is
generated by all $T_b \in {\Bbb K} [{\calS}_r]$ which are 
formed from tensor products {\rm (\ref{equ2.17})} of arbitrary tensors
$T^{(i)} \in {\calT}_{{\frak l}_i^{\ast}}$.
\end{Prop}
Let $\regrep_G : G \rightarrow GL({\Bbb K}[G])$ denote the {\it regular
representation} of a finite group $G$ defined by
$\regrep_g (f) := g \cdot f$, $g \in G$, $f \in {\Bbb K}[G]$. If we use the 
above left ideals ${\frak l}_i$, ${\frak l}$
to define subrepresentations
${\alpha}_i := \regrep_{{\calS}_{r_i}} |_{{\frak l}_i}$,
$\beta := \regrep_{{\calS}_r} |_{\frak l}$,
then the representation $\beta$ is equivalent to a {\itshape Littlewood-Richardson 
product} (see B. Fiedler \cite[Sec.III.3.2]{fie16}):
\begin{eqnarray}
{\frak l} \; {\rm according}\;{\rm to}\;{\rm (\ref{equ2.19})}
& \;\;\Longrightarrow\;\; &
\beta \;\sim\;
{\alpha}_1 \,\#\ldots\#\, {\alpha}_m \uparrow {\calS}_r \label{equ2.20}
\end{eqnarray}
('$\#$' denotes the outer tensor product of the above representations.) This result corresponds to statements of S.A. Fulling et al. \cite{full4}.
Relation (\ref{equ2.20}) allows us to determine information about the structure of the left ideal (\ref{equ2.19}) by means of the {\itshape Littlewood-Richardson rule}\footnote{See the
references \cite{kerber,kerber3,jameskerb,littlew1,mcdonald,full4,fultharr} for the Littlewood-Richardson rule.}.\vspace{20pt}

\section{Proof ot the Theorems \ref{thm1.8}, \ref{thm1.9} and \ref{thm1.10}}
\noindent In this section we prove the Theorems \ref{thm1.8}, \ref{thm1.9} and \ref{thm1.10}.\vspace{10pt}

\subsection{Use of the Littlewood-Richardson rule (Proof of Theorem \ref{thm1.8})} \label{subsec3.1}%
We consider the Problem \ref{probl1.7}. Theorem \ref{thm2.3} and Corollary \ref{cor2.4} tell us that a tensor $W\in\calT_5 V$ is an algebraic covariant derivative curvature tensor iff a tensor $\tilde{W}\in\calT_5 V$ exists such that
$W = \frac{1}{24}\,y_{t'}^{\ast}\tilde{W}$ where $y_{t'}$ is the Young symmetrizer of the tableau (\ref{equ2.16}).

It is well-known that every tensor $\tilde{W}\in\calT_5 V$ can be represented as a finite sum of tensors $W''\otimes W'''$ where $W''\in\calT_2 V$, $W'''\in\calT_3 V$. Thus we can find a finite subset $\calP\subset\calT_2 V\times\calT_3 V$ for every $\tilde{W}\in\calT_5 V$ such that
\begin{eqnarray} \label{equ3.1}
\tilde{W} & = & \sum_{(W'',W''')\in\calP}\,W''\otimes W'''\,.
\end{eqnarray}
Assume that we know defining idempotents $f''\in\bbK[\calS_2]$, $f'''\in\bbK[\calS_3]$ of the symmetry classes of such tensors $W'', W'''$. If we determine decompositions of those idempotents into pairwise orthogonal, primitive idempotents, e.g.
\begin{eqnarray}
f''' & = & \sum_i\,e_i\;\;\;\;\;\;,\;\;\;\;\;\;e_i\cdot e_j\;=\;{\delta}_{i j}\,e_i\;,\;e_i\in\bbK[\calS_3]\;\mathrm{primitive}\,
\end{eqnarray}
then we obtain decompositions of $W''$, $W'''$ into tensors whose symmetry classes are defined by minimal right ideals (e.g. $W''' = f'''W''' = \sum_i\,e_iW'''$). Consequently we can construct such a representation (\ref{equ3.1}) for $\tilde{W}$ in which all $W''$, $W'''$ belong to symmetry classes which are defined by {\itshape minimal right ideals} $\frakr'', \frakr'''$. We determine the structure of the symmetry class of such a product tensor $W''\otimes W'''$.

Because of Proposition \ref{prop2.4} all group ring elements $W''_{b''}, W'''_{b'''}$ according to Definition \ref{def2.3} lie in minimal left ideals $\frakl'' = (\frakr'')^{\ast}$, $\frakl''' = (\frakr''')^{\ast}$. From Proposition \ref{prop2.8} follows that the group ring elements $(W''\otimes W''')_b$ belong to a left ideal $\frakl\subseteq\bbK[\calS_5]$ which is the representation space of a representation $\beta\sim {\alpha}''\,\#\,{\alpha}'''\uparrow\calS_5$ where
${\alpha}'' := \regrep_{\calS_2}|_{\frakl''}$, 
${\alpha}''' := \regrep_{\calS_3}|_{\frakl'''}$. Since $\frakl''$, $\frakl'''$ are minimal the representations ${\alpha}''$, ${\alpha}'''$ are irreducible. Their equivalence classes are characterized by certain partitions, i.e. there are partitions ${\lambda}''\vdash 2$, ${\lambda}'''\vdash 3$, such that
${\alpha}''\sim [{\lambda}'']$,
${\alpha}'''\sim [{\lambda}''']$.

Now we can determine information about the decomposition of $\beta$ into irreducible representations by means of the {\itshape Littlewood-Richardson rule}. There are the following possibilities to form Littlewood-Richardson products:
\begin{eqnarray}
{[3][2]} & \sim & [5] + [3\,2] + [4\,1] \label{equ3.3}\\
{[3][1^2]} & \sim & [4\,1] + [3\,1^2] \nonumber \\
{[2\,1][2]} & \sim & [3\,2] + [4\,1] + [2^2\,1] + [3\,1^2] \label{equ3.4}\\
{[2\,1][1^2]} & \sim & [3\,2] + [2^2\,1] + [3\,1^2] + [2\,1^3] \label{equ3.5}\\
{[1^3][2]} & \sim & [3\,1^2] + [2\,1^3] \nonumber \\
{[1^3][1^2]} & \sim & [2^2\,1] + [2\,1^3] + [1^5] \nonumber
\end{eqnarray}
Every irreducible representation of the right-hand sides of these relations has a representation space which is a minimal left ideal
$\check{\frakl} = \bbK[\calS_5]\cdot\check{f}$ generated by a primitive idempotent $\check{f}\in\bbK[\calS_5]$. A tensor $\check{W}\in\calT_{\check{\frakr}}$ of the symmetry class $\calT_{\check{\frakr}}$ defined by
$\check{\frakr} = {\check{\frakl}}^{\ast}$ satisfies
$\check{W} = {\check{f}}^{\ast}\check{W}$ and the symmetrization of $\check{W}$ by means of $y_{t'}^{\ast}$ yields
$y_{t'}^{\ast}\check{W} = y_{t'}^{\ast}{\check{f}}^{\ast}\check{W} =
(f\cdot y_{t'})^{\ast}\check{W}$. But since the Young symmetrizer $y_{t'}$ generates a minimal left ideal which lies in the equivalence class of minimal left ideals characterized by $(3\,2)\vdash 5$, we obtain $f\cdot y_{t'}\not= 0$ at most if $\check{\frakl}$ belongs to the class of $(3\,2)\vdash 5$, too. This happens only in the cases (\ref{equ3.3}), (\ref{equ3.4}), (\ref{equ3.5}). Thus Theorem \ref{thm1.8} is valid for tensor products (\ref{equ1.18}) of type
$T\otimes\hat{T}$. Tensor products $\hat{T}\otimes T$ can be handled in the same way.\vspace{10pt}

\subsection{Proof of Theorem \ref{thm1.9}} \label{subsec3.2} The proof of Theorem \ref{thm1.9} uses a method which we applied already in the proof of \cite[Theorem 5.1]{fie20}.

If $G\subseteq\calS_r$ is a subgroup of $\calS_r$ and $i_1,\ldots,i_h$ are natural numbers with $1\le i_l\le r$, then we denote by $G_{i_1,\ldots,i_h}$ the {\itshape stabilizer subgroup} of all such permutations from $G$ for which $i_1,\ldots,i_h$ are fixed points.

Now we consider the following group ring elements
\begin{eqnarray}
{\sigma}_k & := & y_{t'}^{\ast}\cdot {\xi}_k\in\bbK[\calS_5]
\;\;\;\;,\;\;\;\;
(k = 1\,,\,2)
\end{eqnarray}
where $y_{t'}$ is the Young symmetrizer of the Young tableau (\ref{equ2.16}) and
\begin{eqnarray}
{\xi}_1 & := & \Big(\sum_{p\in (\calS_5)_{3,4,5}}\,p\Big)\cdot \Big(\sum_{q\in (\calS_5)_{1,2}}\,q\Big)\\
{\xi}_2 & := & \Big(\sum_{q\in (\calS_5)_{4,5}}\,q\Big)\cdot \Big(\sum_{p\in (\calS_5)_{1,2,3}}\,p\Big)
\end{eqnarray}
A calculation (by means of {\sf PERMS} \cite{fie10}) shows that ${\sigma}_1 \not= 0$ and ${\sigma}_2 \not= 0$ (see \cite{fie21}).

The right ideals ${\sigma}_k \cdot \bbK [\calS_5]$ are non-vanishing subideals of the right ideal $\frakr = y_{t'}^{\ast}\cdot \bbK [\calS_5]$ which defines the symmetry class ${\calT}_{\frakr}$ of algebraic covariant derivative curvature tensors. Since $\frakr$ is a minimal right ideal, we obtain
$\frakr = {\sigma}_1 \cdot \bbK [\calS_5] = {\sigma}_2 \cdot \bbK [\calS_5]$, i.e. ${\sigma}_1$ and ${\sigma}_2$ are generating elements of $\frakr$, too.

A tensor $T \in {\calT}_5 V$ is an algebraic covariant derivative curvature tensor iff there exist $a \in \frakr$ and $T' \in {\calT}_5 V$ such that $T = aT'$. Since further every $a \in \frakr$ has representations
$a = {\sigma}_k\cdot x_k$, $k = 1 , 2$, a tensor
$T \in {\calT}_5 V$ is an algebraic covariant derivative curvature tensor iff there is a tensor
$T' \in {\calT}_5 V$ such that $T = {\sigma}_1 T'$ or
$T = {\sigma}_2 T'$.

Let us consider the case $k = 1$. We obtain all algebraic covariant derivative curvature tensors if we form $T = {\sigma}_1 T'$, $T' \in {\calT}_5 V$. As in the proof of Theorem \ref{thm1.8} we can use a representation
\begin{eqnarray}
T' & = & \sum_{(M,N) \in \calP} M \otimes N \label{equ3.9}
\end{eqnarray}
for a $T' \in {\calT}_5 V$, where $\calP$ is a finite set of pairs $(M,N)$ of tensors $M\in\calT_2 V$, $N\in\calT_3 V$ . From \eqref{equ3.9} we obtain
\begin{eqnarray}
{\xi}_1 T' & = & \sum_{(S,\hat{S}) \in {\tilde \calP}} S\otimes\hat{S} \label{equ3.10}
\end{eqnarray}
where ${\tilde \calP}$ is the finite set of pairs $(S,\hat{S})$ of the symmetrized tensors $S = (\sum_{p\in\calS_2}\,p)\,M$, $\hat{S} = (\sum_{q\in\calS_3}\,q)\,N$, $(M,N) \in \calP$.
Now the application of $y_{t'}^{\ast}$ shows, that the algebraic covariant derivative curvature tensor $T$ is a finite sum of tensors 
\[
T \;=\;
{\sigma}_1 T'
\;=\;
(y_t^{\ast}\cdot {\xi}_1) T'
\;=\;
\sum_{(S,\hat{S}) \in {\tilde \calP}} y_{t'}^{\ast}(S\otimes\hat{S})\,.
\]
The case $k = 2$ can be treated in exactly the same way. We can show that
\begin{eqnarray}
{\xi}_2 T' & = & \sum_{(\hat{S},S) \in {\hat \calP}} \hat{S}\otimes S \label{equ3.11}
\end{eqnarray}
where ${\hat \calP}$ is a finite set of pairs $(\hat{S},S)$ of symmetric tensors $S\in\calT_2 V$, $\hat{S}\in\calT_3 V$. Then we obtain
\[
T \;=\;
{\sigma}_2 T'
\;=\;
(y_t^{\ast}\cdot {\xi}_2) T'
\;=\;
\sum_{(\hat{S},S) \in {\hat \calP}} y_{t'}^{\ast}(\hat{S}\otimes S)\,.
\]
Formula (\ref{equ1.21}) can be proved by a computer calculation. We used the {\sf Mathematica} packages {\sf PERMS} \cite{fie10} and {\sf Ricci} \cite{ricci3} to this end (see \cite{fie21}).\vspace{10pt}

\subsection{Proof of Theorem \ref{thm1.10}} \label{subsec3.3}
In the situation considered in Theorem \ref{thm1.10} the symmetry class of the tensors $U$ is not unique. The $(2\,1)$-equivalence class of minimal right ideals $\frakr\subset\bbK[\calS_3]$ which we use to define symmetry classes for the $U$ is an infinite set. The set of generating idempotents for these right ideals $\frakr$ is infinite, too. We use {\itshape discrete Fourier transforms} to determine the family of all primitive idempotents in the $(2\,1)$-equivalence class of minimal right ideals $\frakr\subset\bbK[\calS_3]$.

We denote by $\bbK^{d\times d}$ the set of all $d\times d$-matrices of elements of $\bbK$.
\begin{Def}
A {\it discrete Fourier transform}\footnote{See M. Clausen and U. Baum \cite{clausbaum1,clausbaum2} for details about fast discrete Fourier transforms.} for $\calS_r$ is an isomorphism
\begin{eqnarray}
D : \; \bbK [\calS_r] & \rightarrow &
\bigotimes_{\lambda \vdash r} {\bbK}^{d_{\lambda} \times d_{\lambda}} \label{equ3.12}%
\end{eqnarray}
according to Wedderburn's theorem which maps the group ring $\bbK [\calS_r]$ onto an outer direct product $\bigotimes_{\lambda \vdash r} {\bbK}^{d_{\lambda} \times d_{\lambda}}$ of full matrix rings ${\bbK}^{d_{\lambda} \times d_{\lambda}}$. We denote by $D_{\lambda}$ the {\it natural projections}
$D_{\lambda} : \bbK [\calS_r] \rightarrow
{\bbK}^{d_{\lambda} \times d_{\lambda}}$.
\end{Def}
A discrete Fourier transform maps every group ring element $a\in\bbK[\calS_r]$
to a block diagonal matrix
\begin{eqnarray} \label{equ3.13}%
D :\;\;a\;=\;\sum_{p\in\calS_r}\,a(p)\,p & \mapsto &
\left(
\begin{array}{cccc}
A_{{\lambda}_1} & & & 0 \\
 & A_{{\lambda}_2} & & \\
 & & \ddots & \\
0 & & & A_{{\lambda}_k} \\
\end{array}
\right)\,.
\end{eqnarray}
The matrices $A_{\lambda}$ are numbered by the partitions $\lambda\vdash r$. The dimension $d_{\lambda}$ of every matrix $A_{\lambda}\in{\bbK}^{d_{\lambda} \times d_{\lambda}}$ can be calculated from the Young frame belonging to $\lambda\vdash r$ by means of the {\itshape hook length formula}. For $r = 3$ we have
\begin{center}
\begin{tabular}{c|ccc}
$\lambda$ & $(3)$ & $(2\,1)$ & $(1^3)$ \\
\hline
$d_{\lambda}$ & 1 & 2 & 1 \\
\end{tabular}
\end{center}
The inverse discrete Fourier transform is given by
\begin{Prop}\footnote{See M. Clausen and U. Baum \cite[p.81]{clausbaum1}.}
If $D : \bbK [\calS_r]\rightarrow
\bigotimes_{\lambda \vdash r} {\bbK}^{d_{\lambda} \times d_{\lambda}}$
is a discrete Fourier transform for $\bbK[\calS_r]$, then we have for every
$a\in\bbK[\calS_r]$
\begin{eqnarray}
\forall\,p\in\calS_r:\;\;\;a(p) & = & \frac{1}{r!}\,\sum_{\lambda\vdash r}\,
d_{\lambda}\,\mathrm{trace}\left\{D_{\lambda}(p^{-1})\cdot D_{\lambda}(a)\right\} \label{equ3.14}\\
 & = & \frac{1}{r!}\,\sum_{\lambda\vdash r}\,
d_{\lambda}\,\mathrm{trace}\left\{D_{\lambda}(p^{-1})\cdot A_{\lambda}\right\}\,. \nonumber
\end{eqnarray}
\end{Prop}
It is easy to determine minimal right ideals and primitive idempotents in the image $\bigotimes_{\lambda \vdash r} {\bbK}^{d_{\lambda} \times d_{\lambda}}$ of a discrete Fourier transform. We denote by $\bbK^d$ the vector space of all
$d$-tuples (rows) of elements of $\bbK$.
\begin{Lem}\footnote{See B. Fiedler \cite[p.18,p.20]{fie16}.} We assign to a fixed $d$-tuple $a\in\bbK^d$ the set of $d\times d$-matrices $\frakr_a := \{a^t\cdot u\;|\;u\in\bbK^d\}$. Then we have:
\begin{enumerate}
\item{The sets $\frakr_a$ are minimal right ideals of $\bbK^{d\times d}$ for every $a\in\bbK^d$, $a\not= 0$.}
\item{Let $\frakr\subset\bbK^{d\times d}$ be a minimal right ideal and $A\in\frakr$, $A\not= 0$. Then every non-vanishing column $a^t$ of $A$ yields a $d$-tuple $a\in\bbK^d$ which fulfils $\frakr = \frakr_a$.}
\end{enumerate} \label{lem3.3}%
\end{Lem}
In our considerations we are interested in the matrix ring $\bbK^{2\times 2}$ which corresponds to the $(2\,1)$-equivalence class of minimal right ideals $\frakr\subset\bbK[\calS_3]$. We can formulate the following assertion about generating idempotents of minimal right ideals of $\bbK^{e\times 2}$.
\begin{Prop}
Every minimal right ideal $\frakr\subset\bbK^{2\times 2}$ is generated by exactly one of the following (primitive) idempotents
\begin{eqnarray} \label{equ3.15}%
X_{\nu} & := & \left(
\begin{array}{cc}
1 & 0 \\
\nu & 0 \\
\end{array}
\right)\;\;\;,\;\;\;\nu\in\bbK\,,\\
Y & := & \left(
\begin{array}{cc}
0 & 0 \\
0 & 1 \\
\end{array}
\right)\,. \label{equ3.16}%
\end{eqnarray}
\end{Prop}
\begin{proof}
Lemma \ref{lem3.3} tells us that two right ideals $\frakr_{a_1}$ and $\frakr_{a_2}$ are different iff $a_1$ and $a_2$ are linearly independent. This leads us to the following two cases:
\begin{enumerate}
\item{Every vector $(x,y)\in\bbK^2$ with $x\not= 0$ is proportional to one element of the family
\begin{eqnarray}
(1,\nu)\in\bbK^2 & , & \nu\in\bbK\,. \label{equ3.17}
\end{eqnarray}
Furthermore the elements of the family (\ref{equ3.17}) are pairwise linearly independent.
}
\item{Every vector $(x,y)\in\bbK^2$ with $x = 0$ is proportional to
\begin{eqnarray}
 & & (0,1)\in\bbK^2\,. \label{equ3.18}
\end{eqnarray}
}
\end{enumerate}
Moreover $(0,1)$ and $(1,\nu)$ are linearly independent for all $\nu\in\bbK$. Thus we obtain every minimal right ideal of $\bbK^{2\times 2}$ exactly once if we consider all right ideals $\frakr_a$ where $a$ runs through the family (\ref{equ3.17}) and the element (\ref{equ3.18}). Now we form ansatzes for idempotents of the basis of (\ref{equ3.17}), (\ref{equ3.18}).

\noindent {\bf Case (1):} First we determine a generating idempotent for a minimal right ideal $\frakr_a$ where $a$ belongs to the family (\ref{equ3.17}). We use an ansatz
\begin{eqnarray*}
X_{\nu} & := &
\left(
\begin{array}{c}
1 \\
\nu \\
\end{array}
\right)\cdot (1,y)\;\;\;,\;\;\;\nu\,,\,y\in\bbK\,.
\end{eqnarray*}
Since
$X_{\nu} = \left(
\begin{array}{cc}
1 & y \\
\nu & \nu\,y \\
\end{array}\right)$
and $X_{\nu}\cdot X_{\nu} = \{1 + \nu\,y\}\,X_{\nu}$, we obtain
\[
X_{\nu}\cdot X_{\nu}\;=\;X_{\nu}\;\;\Leftrightarrow\;\;
1 + \nu\,y\;=\;1\;\;\Leftrightarrow\;\;y\;=\;0\,.
\]
In the case $\nu\not= 0$ the only solution of $1 + \nu\,y\;=\;1$ is $y = 0$. In the case $\nu = 0$ values $y\not= 0$ lead to idempotents $X_{\nu}$, too. However, all these idempotents generate one and the same right ideal since the column $(1\;\nu)^t$ is fixed. Because we search for exactly one idempotent in every right ideal, the solution $y = 0$ is sufficient.\\*[0.3cm]
{\bf Case (2):} For an idempotent belonging to (\ref{equ3.18}) we make an ansatz
\begin{eqnarray*}
Y & := &
\left(
\begin{array}{c}
0 \\
1 \\
\end{array}
\right)\cdot (x,y)\;\;\;,\;\;\;x\,,\,y\in\bbK\,.
\end{eqnarray*}
From this ansatz we obtain
$Y = \left(
\begin{array}{cc}
0 & 0 \\
x & y \\
\end{array}\right)$
and $Y\cdot Y = y\,Y$. Consequently we have $y = 1$. Furthermore we can set $x = 0$ because every value of $x$ yields a generating idempotent of one and the same right ideal.
\end{proof}
Now we determine the primitive idempotents of $\bbK[\calS_3]$ which correspond to (\ref{equ3.15}) and (\ref{equ3.16}).
\begin{Prop}
Let us use Young's natural representation\footnote{Three discrete Fourier transforms (\ref{equ3.12}) are known for symmetric groups $\calS_r$: (1) {\itshape Young's natural representation}, (2) {\itshape Young's seminormal representation} and (3) {\itshape Young's orthogonal representation}. (See \cite{boerner,boerner2,kerber,clausbaum1}. A short description of (1) and (2) can be found in \cite[Sec.I.5.2]{fie16}.). All three discrete Fourier transforms are implemented in the program package {\sf SYMMETRICA} \cite{kerbkohn2,kerbkohnlas}. {\sf PERMS} \cite{fie10}  uses the natural representation. The fast {DFT}-algorithm of M. Clausen and U. Baum \cite{clausbaum1,clausbaum2} is based on the seminormal representation.} of $\calS_3$ as discrete Fourier transform. If we apply the Fourier inversion formula {\rm (\ref{equ3.14})} to a $4\times 4$-block matrix
\begin{eqnarray} \label{equ3.19}%
\left(
\begin{array}{ccc}
0 & 0 & 0 \\
0 & A & 0 \\
0 & 0 & 0 \\
\end{array}
\right)
\end{eqnarray}
where $A$ is equal to {\rm (\ref{equ3.15})} or {\rm (\ref{equ3.16})}, then we obtain the following idempotents of $\bbK[\calS_3]$:
\begin{eqnarray}
X_{\nu}\;\;\;\Rightarrow\;\;\;{\xi}_{\nu} & := &
{\textstyle\frac{1}{3}}\,\{[1,2,3] + \nu [1,3,2] + (1-\nu)[2,1,3] \label{equ3.20}\\
 & & - \nu [2,3,1] + (-1+\nu)[3,1,2] - [3,2,1]\} \nonumber \\
Y\;\;\;\Rightarrow\;\;\;\eta & := &
{\textstyle\frac{1}{3}}\,\{[1,2,3] - [2,1,3] - [2,3,1] + [3,2,1]\}\,. \label{equ3.21}
\end{eqnarray}
\end{Prop}
\begin{proof}
We see from formula (\ref{equ3.13}) and the table of the $d_{\lambda}$-values for $\calS_3$ that the images of a discrete Fourier transform of $\calS_3$ are block diagonal matrices
\begin{eqnarray}
D :\;\;a\;=\;\sum_{p\in\calS_r}\,a(p)\,p & \mapsto &
\left(
\begin{array}{ccc}
A_{(3)} & & 0 \\
 & A_{(2\,1)} & \\
0 & & A_{(1^3)} \\
\end{array}
\right)\,,
\end{eqnarray}
where $A_{(3)}$ and $A_{(1^3)}$ are $1\times 1$-matrices and $A_{(2\,1)}$ is a $2\times 2$-matrix. The $(2\,1)$-equivalence class of minimal right ideals of $\bbK[\calS_3]$ corresponds to $A_{(3)} = A_{(1^3)} = 0$. Starting with block matrices (\ref{equ3.19}) and $A = X_{\nu}$ or $A = Y$ we calculated the idempotents ${\xi}_{\nu}$ and $\eta$ by means of the tool
\verb|InvFourierTransform| of the Mathematica package {\sf PERMS} \cite{fie10} (see also \cite[Appendix B]{fie16}.) These calculation can be carried out also by the program package {\sf SYMMETRICA} \cite{kerbkohn2,kerbkohnlas}.
\end{proof}
\begin{Rem}
It is interesting to clear up the connection of the idempotents ${\xi}_{\nu}$ and $\eta$ with Young symmetrizers.
A simple calculation shows that
\begin{eqnarray}
{\xi}_0\;=\;{\textstyle\frac{1}{3}}\,y_{t_1}
 & , &
\eta\;=\;{\textstyle\frac{1}{3}}\,y_{t_2}
\end{eqnarray}
where $y_{t_1}$ and $y_{t_2}$ are the Young symmetrizers of the tableaux
\begin{eqnarray*}
t_1\;=\;
\begin{array}{|c|c|c}
\cline{1-2}
1 & 2 & \\
\cline{1-2}
3 & \multicolumn{2}{c}{\;\;\;} \\
\cline{1-1}
\end{array}
 & , &
t_2\;=\,
\begin{array}{|c|c|c}
\cline{1-2}
1 & 3 & \\
\cline{1-2}
2 & \multicolumn{2}{c}{\;\;\;} \\
\cline{1-1}
\end{array}
\,.
\end{eqnarray*}
\end{Rem}
\noindent {\bf Proof of Theorem \ref{thm1.10}:} The proof ot Theorem \ref{thm1.10} uses the same method as the proof of Theorem \ref{thm1.9}. To treat expressions $y_{t'}^{\ast}(S\otimes U)$ and $y_{t'}^{\ast}(A\otimes U)$ form the following groupring elements of $\bbK[\calS_5]$:
\begin{eqnarray}
{\sigma}_{\nu , \epsilon} & := & y_{t'}^{\ast}\cdot {\zeta}_{\epsilon}'\cdot
{\xi}_{\nu}''
\;\;\;,\;\;\;
{\rho}_{\epsilon}\;:=\;y_{t'}^{\ast}\cdot {\zeta}_{\epsilon}'\cdot
{\eta}'' \label{equ3.24}\\
{\zeta}_{\epsilon}' & := & \id + \epsilon\,(1\,2)\;\;\;,\;\;\;
\epsilon\in\{1,-1\}\\
{\xi}_{\nu} & \mapsto & {\xi}_{\nu}''\in\bbK[\calS_5]
\;\;\;,\;\;\;
{\eta}\;\mapsto\;{\eta}''\in\bbK[\calS_5]\,. \label{equ3.26}
\end{eqnarray}
Formula (\ref{equ3.26}) denotes the embedding of the group ring elements ${\xi}_{\nu}\,,\,\eta\in\bbK[\calS_3]$ into $\bbK[\calS_5]$ which is induced by the mapping
\begin{eqnarray*}
\calS_3\rightarrow\calS_5 & , & 
[i_1,i_2,i_3]\mapsto [1,2,i_1 + 2,i_2 + 2,i_3 + 2]\,.
\end{eqnarray*}
For expressions $y_{t'}^{\ast}(U\otimes S)$ and $y_{t'}^{\ast}(U\otimes A)$ we consider the groupring elements of $\bbK[\calS_5]$:
\begin{eqnarray}
{\sigma}_{\nu , \epsilon} & := & y_{t'}^{\ast}\cdot {\zeta}_{\epsilon}''\cdot
{\xi}_{\nu}'
\;\;\;,\;\;\;
{\rho}_{\epsilon}\;:=\;y_{t'}^{\ast}\cdot {\zeta}_{\epsilon}''\cdot
{\eta}' \label{equ3.27}\\
{\zeta}_{\epsilon}'' & := & \id + \epsilon\,(4\,5)\;\;\;,\;\;\;
\epsilon\in\{1,-1\}\\
{\xi}_{\nu} & \mapsto & {\xi}_{\nu}'\in\bbK[\calS_5]
\;\;\;,\;\;\;
{\eta}\;\mapsto\;{\eta}'\in\bbK[\calS_5]\,. \label{equ3.29}
\end{eqnarray}
Formula (\ref{equ3.29}) denotes the embedding of the group ring elements ${\xi}_{\nu}\,,\,\eta\in\bbK[\calS_3]$ into $\bbK[\calS_5]$ which is induced by the mapping
\begin{eqnarray*}
\calS_3\rightarrow\calS_5 & , & 
[i_1,i_2,i_3]\mapsto [i_1,i_2,i_3,4,5]\,.
\end{eqnarray*}
Using the Mathematica package {\sf PERMS} \cite{fie10} we verified 
\begin{center}
\fbox{
${\rho}_{\epsilon}\not= 0 \;\;\;\;\;\;\mathrm{and}\;\;\;\;\;\;
{\sigma}_{\nu , \epsilon}\not= 0\;\Leftrightarrow\;\nu\not=\frac{1}{2}
$}
\end{center}
in both cases. Thus the same arguments which we used in the proof of Theorem \ref{thm1.9} yield now in both cases:
\begin{enumerate}
\item{The minimal right ideals $y_{t'}^{\ast}\cdot\bbK[\calS_5]$, ${\rho}_{\epsilon}\cdot\bbK[\calS_5]$ and  ${\sigma}_{\nu , \epsilon}\cdot\bbK[\calS_5]$
coincide, if $\nu\not=\frac{1}{2}$.}
\item{A tensor $T\in\calT_5 V$ is an algebraic covariant derivative curvature tensor iff a tensor $T'\in\calT_5 V$ exists such that
$T = {\rho}_{\epsilon} T'$ or $T = {\sigma}_{\nu , \epsilon} T'$ (if $\nu\not=\frac{1}{2}$).
}
\end{enumerate}
Now we represent the tensor $T'$ by finite sums of product tensors. We use sums
\begin{eqnarray}
T' & = & \sum_{(M,N)\in\calP}\,M\otimes N\;\;\;,\;\;\;
\calP\subset\calT_2 V\times\calT_3 V\;\mathrm{finite}
\end{eqnarray}
in the case (\ref{equ3.24})--(\ref{equ3.26}) and sums
\begin{eqnarray}
T' & = & \sum_{(N,M)\in\calP}\,N\otimes M\;\;\;,\;\;\;
\calP\subset\calT_3 V\times\calT_2 V\;\mathrm{finite}
\end{eqnarray}
in the case (\ref{equ3.27})--(\ref{equ3.29}).

Now the argumentation of the proof of Theorem \ref{thm1.9} yield also the statement of Theorem \ref{thm1.10} (if $\nu\not=\frac{1}{2}$). The symmetrizations
${\zeta}_{\epsilon}'\cdot {\xi}_{\nu}''(M\otimes N)$,
${\zeta}_{\epsilon}'\cdot {\eta}''(M\otimes N)$, 
${\xi}_{\nu}'\cdot {\zeta}_{\epsilon}''(N\otimes M)$ and
${\eta}'\cdot {\zeta}_{\epsilon}''(N\otimes M)$ lead to product tensors
$S\otimes U$, $A\otimes U$, $U\otimes S$ and $U\otimes A$.
Note that the value $\epsilon = 1$ ($\epsilon = -1$) produces symmetric tensors $S\in\calT_2 V$ (alternating tensors $A\in\calT_2 V$) from the tensors $M\in\calT_2 V$.

Thus we can write every tensor
${\zeta}_{\epsilon}'\cdot {\xi}_{\nu}''T'$,
${\zeta}_{\epsilon}'\cdot {\eta}''T'$, 
${\xi}_{\nu}'\cdot {\zeta}_{\epsilon}''T'$ and
${\eta}'\cdot {\zeta}_{\epsilon}''T'$
as a finite sum of suitable tensors 
$S\otimes U$, $A\otimes U$, $U\otimes S$ or $U\otimes A$.
An application of $y_{t'}^{\ast}$ shows that every algebraic covariant derivative curvature tensor $T$ can be expressed by a finite sum of each of the tensor types (\ref{equ1.22})

Finally we can verify by a simple calculation that the idempotent ${\xi}_{\nu}$ is equal to the idempotent (\ref{equ1.23}) if $\nu = \frac{1}{2}$. This finishes the proof of Theorem \ref{thm1.10}.

\vspace{1cm}

\noindent {\bf Acknowledgements.} I would like to thank Prof. P. B. Gilkey for important and helpful discussions and for valuable suggestions for future investigations.
\vspace{0.4cm}

\end{document}